\documentclass{article}
\usepackage{amsfonts}
\usepackage{amssymb}
\usepackage{latexsym}
\usepackage{ifthen}

\newtheorem{thm}{Theorem}
\newtheorem{lem}{Lemma}

\newtheorem{cor}{Corollary}
\newcommand{\ra}{\rightarrow}
\newcommand{\la}{\leftarrow}

\newcommand{\da}{\downarrow}
\newcommand{\tx}{\textrm}

\begin{document}

\title{KW-models for (multivariate) linear differential systems}
\date{~~}
\maketitle
\author\begin{center}{Vakhtang Lomadze}\end{center}
\begin{center}\footnotesize{\textsl{ {
Department of Mathematics, I. Javakhishvili Tbilisi State
University, Tbilisi 0183, Georgia }}}
\end{center}

{\bf Abstract}. A class of state models, called Kronecker-Weierstrass models (or, simply, KW-models), is introduced, and the state representation problem for  linear differential systems is studied in the context of these models.  It is shown, in particular, that
there is a canonical one-to-one correspondence between linear differential systems and similarity classes of minimal proper KW-models.

\bigskip


\section{Introduction}

Throughout, $\mathbb{F}$ is the  field of real or complex numbers,  $s=(s_1,\ldots ,s_n)$ is a sequence of indeterminates and $s_0$ an extra inderterminate,
 $q$ is a fixed positive integer.
We let   ${\mathcal U}$ be the space of infinitely differentiable functions (or distributions) and
$\partial=(\partial_1,\ldots ,\partial_n)$
 the sequence of  partial differential operators.

By a linear differential system (with signal number $q$), one  understands a subset of ${\mathcal U}^q$ that can
be represented as the solution set of an equation of the form
$$
R(\partial)w=0 \ \ \ (w\in {\mathcal U}^q),
$$
where $R$ is a polynomial matrix with $q$ columns. Every such a polynomial matrix is called an AR-representation.

 In this paper, we are concerned with the question of  finding a state space model of {\em minimal}
size implementing a given linear differential system.

In the case when $n=1$, the question is well-studied and understood. We recall (see, for example,  Willems \cite{W}) that every linear differential system can be represented
by a  first order  ordinary differential equation of the form
\begin{equation}
Kx^\prime-Lx=Mw,
\end{equation}
where $K, L, M$ are scalar matrices such that

(a) $\lambda_1K-\lambda_0L$ has full column rank for all $0\neq (\lambda_0,\lambda_1)\in \mathbb{C}^2$;

(b) $[\lambda_1K-\lambda_0L \ | \ M]$  has full row rank when $(\lambda_0,\lambda_1)=(0,1)$.\\
Such a representation is unique  up to isomorphism, and its
 size is as small as possible.

 As is known, condition (a) characterizes  observability, and condition (b) corresponds to "controllability at infinity".
 (We remind that controllability is  characterized by the condition that $[\lambda_1K-\lambda_0L \ | \ M]$  is of  full row rank for all $0\neq (\lambda_0,\lambda_1)\in \mathbb{C}^2$
  (see Proposition VII.11 in Willems \cite{W}).)

By the classical Kronecker-Weierstrass theorem, any pencil $sK-L$ that satisfies  condition (a)
 is similar to a direct sum of Jordan pencils and
{\em vice versa}. (The Jordan pencil of size $d$ is the pencil
$$
s_1\left[\begin{array}{ccc}
 1&  &  \\
  &\ddots &  \\
 &   & 1\\
 0& \cdots  & 0
\end{array}\right] - \left[\begin{array}{ccc}
 0&\cdots  & 0 \\
 1 & &  \\
 & \ddots  & \\
 &   & 1
\end{array}\right],
$$
where the scalar matrices have size $(d+1)\times d$.)

As a natural $n$-dimensional analog of (1), we consider
 a first order partial differential equation
\begin{equation}
K_1\partial_1(x)+\cdots + K_n\partial_n(x)-Lx=Mw
\end{equation}
that satisfies   the following two conditions:

(A) the pencil $s_1K_1+\cdots + s_nK_n-L$ is a direct sum of Jordan pencils;

(B) the matrix $\left[\begin{array}{cc}s_1K_1+\cdots + s_nK_n-s_0L & M\end{array}\right]$ is  regular at infinity.\\
 (The variables $x$ and $w$ above are respectively the state  and  external variables.)

An equation of the form (2) satisfying condition (A) is called a KW-model.
According to \cite{L3}, linear differential systems can always be represented by such models, and this motivates their study.
(The definition of Jordan pencils in $n$ variables will be given later, in Section 3.)

The  condition (B)  should be interpreted as the condition of "controllability at infinity".
As in \cite{L2}, here also we argue  that while representing a linear differential system infinity is needed to be taken into account.
Linear differential systems {\em a priori} are "controllable at infinity", and therefore we only  allow  representations that
are regular at infinity.
KW-models satisfying condition (B) are said to be {\em proper}.

In this paper, we show that every linear differential system has a proper KW-representation and that all {\em minimal}
proper KW-representations are isomorphic. (As in the classical case, minimality means having the size as small as possible.)

{\em Remark}. In dimension 1, Jordan pencils have full column rank, and that is why proper KW-models automatically are minimal. This is not true in higher dimensions.

We put
$S=\mathbb{F}[s]$ and $T=\mathbb{F}[s_0,s]$.
For every  $d\in \mathbb{Z}$, we  write $S_{\leq d}$ to denote the space of polynomials of degree
$\leq d$ and $T_d$ for  the space of homogeneous polynomials of degree
$d$. (For negative $d$,  $S_{\leq d}=\{0\}$ and $T_d=\{0\}$.) Remark that for every $d$, there is a canonical isomorphism
\begin{equation}
S_{\leq d}\simeq T_d.
\end{equation}

For a positive integer $p$,  $[1,p]$  stands for the set
$\{1,\ldots ,p\}$.
If $k$ is a nonnegative  integer, then $\Delta(k)$ denotes the set of multi-indices $i\in \mathbb{Z}_+^n$
of order less than or equal to $k$. (The order of $i=(i_1,\ldots , i_n)$ is defined  to be $|i|=i_1+\cdots + i_n$.)
 We let  $\tau_1, \ldots ,\tau_p$ denote the  partial forward shift operators in the space of all compactly supported $\mathbb{F}$-valued functions on $\mathbb{Z}_+^n$.
 There is an obvious  natural isomorphism of the mentioned space onto $S$, which is given by
 $$
 a\mapsto \sum a(i_1,\ldots ,i_n)s_1^{i_1}\ldots s_n^{i_n};
 $$
 we shall denote it by $\phi$. Clearly, for every $k\geq 0$, $$\phi(\mathbb{F}^{\Delta(k)})=S_{\leq k}.$$
(We identify functions in $\mathbb{F}^{\Delta(k)}$ with functions that are defined on $\mathbb{Z}_+^n$ and have
 support in $\Delta(k)$.)
Remark  that
$$
\phi\tau_k=s_k\phi
$$
 for each $k=1,\ldots , n$.

\section{Preliminaries}

A graded  module over  $T$ is a module $M$ together with a gradation, i.e., a collection $M_d, d\in \mathbb{Z}$ of $\mathbb{F}$-linear subspaces of $M$
such that
$$
M=\bigoplus_{d\in \mathbb{Z}}M_d \ \ \ \tx{and}\ \ \ s_kM_d\subseteq M_{d+1}\ \forall k,d.
$$
The elements of $M_d$ are called the homogeneous elements of $M$ of degree $d$.
A submodule $N\subseteq M$ is called a graded submodule of $M$ if $N=\bigoplus (N\cap M_d)$. Notice that a submodule
 is graded if it contains all the homogeneous components
of each of its elements.

For a graded $T$-module $M$ and an integer
$k$, one denotes by $M(k)$ the graded $T$-module whose
 homogeneous components are defined by
$$M(k)_d = M_{k+d}.$$

{\em Example 1}.  Let $l$ be a positive integer. Then, a sequence $a=(a_1,\ldots ,a_l)$ of integers
determines on $T^l$  a  gradation  consisting of the spaces
$$
T^l(a)_d=\{f\in T^l |\ \deg(f_i)= d+a_i \} \ \ \ (d\geq 0).
$$
The module $T^l$ equipped with this gradation is denoted by
$T^l(a)$. Notice that $T^l(a)=\oplus T(a_i)$.

If $M$ is a graded module, one denotes
by  $M_+$  the graded
 module whose component in the degree $d$ is $M_d$ if $d\geq 0$ and is zero otherwise.

A homogeneous homomorphism
$\varphi :M \ra N$ of graded $T$-modules of degree $d$ is a $T$-module homomorphism
such that $\varphi(M_i)\subseteq N_{i+d}$ for all $i$. The space of all homogeneous homomorphisms of degree $d$
is denoted by $Hom^d(M,N)$.

Homogeneous homomorphisms of degree 0
simply are called homogeneous homomorphisms. These are homomorphisms that preserve the grading.

If $M, N$ are graded modules and $\varphi: N\ra N$ is a homogeneous homomorphism,
then $Ker(\varphi)$ is a graded submodule of $M$ and $Im(\varphi)$ is a graded submodule of $N$.
Moreover, $Coker(\varphi)$ is naturally graded with components $N_d/(Im(\varphi)_d$. (Graded modules together with
homogeneous homomorphisms form an abelian category.)

{\em Example 2}. Let $a$ and $b$ be nonnegative integers. Homogeneous homomorphisms from $T(a)$ to $T(b)$ are exactly multiplications by homogeneous polynomials of degree $b-a$. That is,
$$Hom^0(T(a),T(b))=T_{b-a}.$$

 A polynomial matrix with entries in $T$ is called row-homogeneous (resp. column-homogeneous)  if all the entries in each row
  (resp. column) are homogeneous and have the same degree.

{\em Example 3}. i) A row-homogeneous polynomial matrix $H$ of size $p\times q$ and with row degree  $a=(a_1,\ldots ,a_p)$
determines a homogeneous homomorphism of graded modules
$$
H: T^q\ra T^p(a).
$$
ii)  A column-homogeneous polynomial matrix $H$ of size $q\times p$ and with column degree $a=(a_1,\ldots ,a_q)$
determines a homogeneous homomorphism of graded modules
$$
H: T^p(-a)\ra T^q.
$$

One knows that if $M$ is a finitely generated graded module, then $Hom(M,N)$ has a natural grading for every graded module $N$. Its homogeneous component of degree
$d$ consists of homogeneous homomorphisms of degree $d$, that is,
$$
Hom(M,N)_d=Hom^d(M,N).
$$
 In particular, the dual module $M^*=Hom(M,T)$ is a graded module.

{\em Example 4}. For every integer $k$,
$$
T(k)^*=T(-k).
$$

{\em Example 5}. If $M$ is a finitely generated graded module with $M_{-d}$=\{0\} for all sufficiently large $d$, then
$$
(M_+)^*=M^*.
$$
(This is because $M_+=M/(\oplus_{d<0} M_d)$ and $\oplus_{d<0} M_d$ is a torsion module.)

A graded module that is free as a $T$-module with a basis of homogeneous
elements is called a free graded module. Every (nonzero) free graded $T$-module is isomorphic to a module of the form $T^l(a)$, where
 $l$ is a positive integer and $a$ a sequence of integers of length $l$. Such a representation of a free graded module is unique as the following example says.

 {\em Example 6}. If
$l$, $m$ are positive integers and $a$, $b$ are sequences of  integers respectively of lengths $l$, $m$, then
$$
T^l(a)\simeq T^m(b) \ \ \Rightarrow\ \ l=m\ \ \tx{and}\ \ a=b \ \tx{(up \ to\ permutation)}.
$$
(We do not know a reference for this well-known fact. A proof can be adapted from the proof of Lemma 1 in \cite{L3}.)

{\em Remark}. Everything above is extended with no changes to graded modules over $S$.   It is worth noting that by tensoring a graded $T$-module with $S=T/s_0T$ one gets
a graded $S$-module. In particular, for every integer $d$, we have $T(d)\otimes S=S(d)$. We remark also that if $Q$ is a matrix with entries in $T$, then
$Q\otimes S=Q(0,s)$.

If $R$ is a polynomial matrix of size $p\times q$  with entries in $S$, then one
defines
the homogenization $R^h$ of $R$  by setting $$R^h=diag(s_0^{a(1)}, \ldots , s_0^{a(p)})R(s_1/s_0,\ldots ,s_n/s_0),$$ where $a_1,\ldots ,a_p$ are the row degrees of $R$. Clearly, $R^h$ is a row-homogeneous polynomial matrix (with entries in $T$).
Conversely, if $Q$ is a row-homogeneous (or column-homogeneous) polynomial matrix with entries in $T$, one defines the dehomogenization $Q^{dh}$ by setting
$Q^{dh}=Q(1,s)$.

\section{Kronecker-Weierstrass pencils}

A pencil is a diagram
$$ X\stackrel{K,L}{ \longrightarrow} Y,$$
 where $X, Y$ are (finite-dimensional) linear spaces,
$K=(K_1, \ldots ,K_n)$ is an $n$-tuple of linear maps from $X$ to $Y$  and $L$ is a linear map from $X$ to $Y$ as well.  It gives rise to a  homomorphism
$$
sK -L: X\otimes S\ra  Y\otimes S
$$
given by
$$x\otimes f \mapsto \ K_1(x)\otimes s_1f +\cdots +  K_n(x)\otimes s_nf-L(x)\otimes f.$$

 Two pencils $X\stackrel{K,L}{ \longrightarrow} Y$ and $ X^\prime\stackrel{K^\prime,L^\prime}{ \longrightarrow} Y^\prime$ are
said to be isomorphic (or similar) when there exist isomorphisms $V:X\ra X^\prime$ and $U:Y\ra Y^\prime$
 such that
$$
(sK^\prime -L^\prime)V=U(sK -L).
$$

The simplest (non-trivial) examples of pencils are Jordan pencils.

Let $d$ be a nonnegative integer. A Jordan pencil of index $d$ is the pencil
$$ X^{(d)}\stackrel{I^{(d)},J^{(d)}}{ \longrightarrow} \mathbb{F}^{\Delta(d)},$$
where
$$
X^{(d)}=\{\left(\begin{array}{c}
x_1\\\vdots \\x_n
\end{array}\right)\in (\mathbb{F}^{\Delta(d)})^n \ |\  \tau_1(x_1)+\cdots +\tau_n(x_n) \in \mathbb{F}^{\Delta(d)}\},
$$
and   $I^{(d)}_i:X^{(d)}\ra \mathbb{F}^{\Delta(d)}$ is a linear map defined  by
$$
\left(\begin{array}{c}
x_1\\\vdots \\x_n
\end{array}\right)\mapsto x_i \ \ \ \ \ (i=1,\ldots , n),
$$
and  $J^{(d)}: X^{(d)}\ra \mathbb{F}^{\Delta(d)}$ a linear map defined by
$$
\left(\begin{array}{c}
x_1\\\vdots \\x_n
\end{array}\right)\mapsto \tau_1(x_1)+\cdots +\tau_n(x_n).
$$

\begin{lem} We have
$$
\dim X^{(d)}= \frac{nd}{d+1}\cdot \left(\begin{array}{c} n+d  \\
d
\end{array}\right) \ \ \ \tx{and}\ \ \  \dim \mathbb{F}^{\Delta(d)}= \left(\begin{array}{c}
 n+d  \\
d
\end{array}\right).
$$
\end{lem}
{\em Proof}. The second formula is well-known. To derive the first one, consider the linear map  $(S_{\leq d})^{n+1}\ra S_{d+1}$  given by
$$\left(\begin{array}{c}f_0\\
f_1\\\vdots \\f_n
\end{array}\right)\mapsto s_1f_1+\cdots +s_nf_n-f_0.
$$
This clearly is surjective and its kernel is isomorphic to $ X^{(d)}$;
more precisely, we have the following canonical exact sequence
 \begin{equation}
0\ra X^{(d)}\stackrel{\left[\begin{array}{l}\phi L\\
\phi K_1\\ \ \ \vdots \\ \phi K_n
\end{array}\right]}{\longrightarrow} (S_{\leq d})^{n+1} \stackrel{\left[\begin{array}{cccc}-1& s_1 &\ldots &s_n\end{array}\right]}
{\longrightarrow} S_{\leq d+1}\ra 0.
 \end{equation}

The formula follows.  $\quad\Box$

\begin{lem} There is a canonical exact sequence
$$
0\ra T(-d)\ra  \mathbb{F}^{\Delta(d)}\otimes T\ra  (X^{(d)})^*\otimes T(1).
$$
\end{lem}
{\em Proof}.
We have proved in \cite{L3} that, for each integer $k\geq 0$, the sequence
$$
 X^{(d)}\otimes S_{\leq k-1}\ra \mathbb{F}^{\Delta(d)}\otimes S_{\leq k}\ra S_{\leq d+k}\ra 0
$$
is exact.
In view of (3), it follows from this that the sequence
$$
 X^{(d)}\otimes T_{k-1}\ra \mathbb{F}^{\Delta(d)}\otimes T_k\ra T_{d+k}\ra 0
$$
is exact  for every $k\geq 0$. That is, the sequence
 of graded modules
\begin{equation}
 X^{(d)}\otimes T(-1)\ra \mathbb{F}^{\Delta(d)}\otimes T\ra T(d)_+\ra 0
 \end{equation}
is exact.
By Examples 4 and 5, $ (T(d)_+)^*=T(-d)$, and dualizing we obtain  the required exact sequence.

 The proof is complete. $\quad\Box$

\begin{cor} There is a canonical exact sequence
$$
0\ra S(-d)\ra  \mathbb{F}^{\Delta(d)}\otimes S\ra  (X^{(d)})^*\otimes S(1).
$$
\end{cor}
{\em Proof}. Tensoring (5), by $S=T/s_0T$, we get an exact sequence
$$
 X^{(d)}\otimes S(-1)\ra \mathbb{F}^{\Delta(d)}\otimes S\ra S(d)_+\ra 0.
$$
Dualizing this, as above, we get what we want. $\quad\Box$

One defines in an obvious way direct sums of pencils. Given a sequence $d=(d_1,\ldots ,d_p)$ of integers, let
$KW(d)$ denote the direct sum of Jordan pencils of sizes $d_1,\ldots ,d_p$.

One consequence of the previous lemma  is the following lemma.
\begin{lem} Let $a=(a_1,\ldots ,a_l)$ and $b=(b_1,\ldots ,b_m)$ be  two sequences of nonnegative integers. Then
$$
KW(a)\simeq KW(b)\ \ \Rightarrow\ \ l=m\ \ \tx{and}\ \ a=b \   \ (\tx{up\ to\ permutation}).
$$
\end{lem}
{\em Proof}. Using Lemma 2, we can see that the isomorphism $KW(a)\simeq KW(b)$ yields an isomorphism
$T^l(-a)\simeq T^m(-b)$.

The lemma follows (see Example 6).  $\quad\Box$

 Call
  a KW-pencil any pencil that  is similar to a direct sum of Jordan pencils.
   By the above lemma,
 the decomposition of a KW-pencil into a direct sum of Jordan pencils  is unique (up to permutation). The sizes of these Jordan pencils are called
the Kronecker indices.

By Lemma 1, if $X\stackrel{K,L}{ \longrightarrow} Y$ is a KW-pencil with Kronecker indices $d_1,\ldots ,d_p$, then
\begin{equation}
\dim(X)=\sum_i \frac{nd_i}{d_i+1}\cdot \left(\begin{array}{c}
 n+d_i  \\
d_i
\end{array}\right)
 \ \  \tx{and}\ \
 \dim(Y)=\sum_i\left(\begin{array}{c}
 n+d_i  \\
d_i
\end{array}\right).
\end{equation}

\section{KW-models and their behaviors}

 A KW-model (with signal number $q$) is a pair consisting of a KW-pencil
 $X\stackrel{K,L}{\longrightarrow} Y$  and  a linear map  $M: \mathbb{F}^q\ra Y$. We can  express it by means of
the diagram
$$
 X\stackrel{K,L}{\longrightarrow} Y \stackrel{M}{\la} \mathbb{F}^q.
 $$
We call $X$ the space of  state variables and $Y$ the space of auxiliary variables. One defines in an obvious way the Kronecker indices of a KW-model.

 Two KW-models $X\stackrel{K,L}{ \longrightarrow} Y\stackrel{M}{\la} \mathbb{F}^q$ and $ X^\prime\stackrel{K^\prime,L^\prime}{ \longrightarrow} Y^\prime
 \stackrel{M^\prime}{\la} \mathbb{F}^q$ are
said to be isomorphic (or similar) when there exist isomorphisms $V:X\ra X^\prime$ and $U:Y\ra Y^\prime$
 such that
$$
(sK^\prime -L^\prime)V=U(sK -L)\ \ \tx{and}\ \ M^\prime= UM.
$$

KW-models appear naturally.
As shown in  \cite{L3},  associated with every AR-model $R$ there is a  canonical  KW-model, namely,
$$
KW(R):\ \ \  \ (KW(a),H^0(R)),
$$
where $a=(a_1,\ldots ,a_p)$ is the sequence of row degrees of $R$ and $H^0(R)$ is the composition $$\mathbb{F}^q\stackrel{R}{\ra}\oplus S_{\leq a(i)}\stackrel{\phi^{-1}}{\ra}\oplus \mathbb{F}^{\Delta(a(i))}.$$
({\em Remark}. The linear map $H^0(R)$ has a cohomological interpretation, and this justifies the notation.)

Assume we have a KW-model
$$
\Sigma:\ \ \  X\stackrel{K,L}{\longrightarrow} Y \stackrel{M}{\la} \mathbb{F}^q.
$$
Tensoring
$$
 X\otimes S\stackrel{sK - L}{\longrightarrow} Y\otimes S \stackrel{M}{\la} S^q
 $$
  with ${\cal U}$ yields the diagram
$$
 X\otimes {\cal U}\stackrel{K\partial - L}{\longrightarrow} Y\otimes {\cal U} \stackrel{M}{\la} {\cal U}^q.
 $$
 (By $K\partial-L$, we mean the differential operator
$ X\otimes {\cal U}\ra  Y\otimes {\cal U}
$
defined by
$$x\otimes u \mapsto \ K_1(x)\otimes \partial_1(u) +\cdots +  K_n(x)\otimes \partial_n(u)-L(x)\otimes u.)$$
We thus have  a partial differential equation
$$
K\partial(x)-L x=Mw\ \ \ (x\in X\otimes {\cal U}, \ w\in {\cal U}^q).
$$
Following Willems \cite{W}, we define the manifest behavior of $\Sigma$ by setting
$$
Bh(\Sigma)=
\{w\in {\cal U}^q |\ \exists x\in X\otimes {\cal U} \ \tx{such \ that} \ K\partial(x)-L x=Mw\}.
$$

The following theorem should be compared with Theorem VII.4 in the a.m. paper by Willems.

\begin{thm}  The manifest behavior of a KW-model is  a linear differential system   and conversely, every linear differential system can
be represented as the manifest behavior of a KW-model.
\end{thm}
{\em Proof}. Let $\Sigma$ be as above. By Lemma 4 in \cite{L3}, there is an exact sequence
$$
 X\otimes S\ra Y\otimes S\ra S^p\ra 0.
$$
Let $R$ denote the composition $S^q \stackrel{M}{\ra} Y\otimes S \ra S^p$. We then have  a commutative diagram
$$
\begin{array}{ccccc}
&  S^q & = & S^q & \\
& M \da\ \ \ \ &  &\ \ \ \da R& \\
X\otimes S\ra  & Y\otimes S & \ra & S^p& \ra 0
\end{array}.
$$
Tensoring this diagram with $\cal U$, we get the following commutative diagram
$$
\begin{array}{ccccc}
&  {\cal U}^q & = &{\cal U}^q & \\
&M \da\ \ \ \ \ &  &\ \ \ \ \ \ \da R(\partial) & \\
X\otimes {\cal U}\ra  & Y\otimes {\cal U} & \ra & {\cal U}^p & \ra 0
\end{array}.
$$
The bottom row here is exact, and this implies  that the manifest behavior of our model is $KerR(\partial)$.

The converse statement is immediate from Theorem 1 in \cite{L3}.

The proof is complete. $\quad\Box$

{\em Remark}. For a KW-model $\Sigma$, we have constructed in the proof of the lemma an AR-model $R$ such that $Bh(\Sigma)=Bh(R)$.
This, however, does not mean that $\Sigma=KW(R)$. The point is that the row degrees
of $R$ are less than or equal to the Kronecker indices of $\Sigma$, but not equal necessarily.

\section{Properness}

 Let $\varphi: F\ra G$ be a homogeneous homomorphism of free graded  modules (over $T$). If   $s_0$ is not zero divisor on the quotient module
$$
F^*/\varphi^*(G^*),
$$
we say that $\varphi$ is {\em regular at infinity}.

The following lemma is a graded (abstract) version of the PBH test for (multivariate)  linear differential systems.

\begin{lem} Let
$$
0\ra F_l \ra \cdots \ra  F_2 \ra  F_1  \stackrel{\varphi}{\ra}  F_0
$$
be an exact sequence of torsion-free graded $T$-modules (and homogeneous homomorphisms), and let $\lambda$ be a homogeneous polynomial.  Then, $\lambda$ is not zero divisor on the quotient module
$F_0/Im(\varphi)$ if and only if the sequence
$$
0\ra F_l/\lambda F_l \ra \cdots \ra  F_2/\lambda F_2 \ra  F_1/\lambda F_1  \ra  F_0/\lambda F_0
$$
is exact.
\end{lem}
{\em Proof}. The proof is the same as that of Corollary 1 in \cite{L1}. $\quad\Box$

  Let a KW-model  $X\stackrel{K,L}{\longrightarrow} Y \stackrel{M}{\la} \mathbb{F}^q$ be given. We say that it is  proper  if the homogeneous homomorphism
$$
\left[\begin{array}{cc}sK-s_0L & M\end{array}\right]: T(-1)\otimes X\oplus T^q \ra T\otimes Y
$$
is regular at infinity.

Let $d_1,\ldots ,d_p$ be the Kronecker indices of our model, and  let $\widetilde{M}$ denote the composition $T^q \stackrel{M}{\ra} Y\otimes T \ra \oplus T(d_k)$. We then  have  a commutative diagram
$$
\begin{array}{ccc}
  T^q & \ra & Y\otimes T \\
 || &  & \da \\
  T^q & \ra & \oplus T(d_k)
\end{array}.
$$
Dualizing this diagram, we get the following commutative diagram
$$
\begin{array}{ccc}
 \oplus T(-d_k) & \ra &  T^q \\
 \da &  & || \\
  Y^*\otimes T & \ra & T^q
\end{array}.
$$
Next, by Lemma 2, we have an exact sequence
$$
0\ra \oplus T(-d_k)\ra T\otimes Y^*\stackrel{sK^t -s_0L^t}{\longrightarrow} T(1)\otimes X^*.
$$
We can easily see that the two homomorphisms
$$
\left[\begin{array}{c}sK^{tr}-s_0L^{tr}\\M^{tr}\end{array}\right]: T\otimes Y^* \ra T(1)\otimes X^* \oplus T^q\ \ \ \tx{and}\ \ \ \widetilde{M}^{tr}: \oplus T(-d_k) \ra  T^q
$$
have isomorphic kernels. Letting $E_1$ and $E_2$ denote these kernels, we have
 a commutative diagram
\begin{equation}
\begin{array}{cccccc}
0\ra & E_1 & \ra & \oplus T(-d_k) & \ra & T^q\\
     & \da &     &     \da        &     & \da   \\
0\ra & E_2 & \ra & T\otimes Y^*   & \ra & T(1)\otimes X^* \oplus T^q
\end{array},
\end{equation}
where the rows are exact and the left downward arrow is an isomorphism.

\begin{lem} $\Sigma$ is proper if and only if \ $T^q \stackrel{\widetilde{M}}{\ra}  \oplus T(d_k)$ is regular at infinity.
\end{lem}

{\em Proof}. By the above lemma, $\Sigma$ is proper if and only if the complex
$$
0 \ra E_1/s_0E_1  \ra S\otimes Y^*   \stackrel{\left[\begin{array}{c}sK^{tr}\\M^{tr}\end{array}\right]}{\ra}  S(1)\otimes X^* \oplus S^q
$$
is exact. Similarly, $T^q \stackrel{\widetilde{M}}{\ra}  \oplus T(-d_k)$ is regular at infinity if and only if the complex
$$
0 \ra E_2/s_0E_2 \ra  \oplus S(-d_k) \stackrel{\widetilde{M}^{tr}(0,s)}{\ra}  S^q
$$
is exact.

We claim that the above two complexes  have the same homologies. Indeed, tensoring (7) by $S=T/s_0T$ yields
 the following commutative diagram
$$
\begin{array}{cccccc}
0 \ra &E_1/s_0E_1 & \ra & \oplus S(-d_k) & \ra & S^q\\
      &     \da      &     &     \da        &                  & \da   \\
0 \ra &E_2/s_0E_2 & \ra & S\otimes Y^*   & \ra &  S(1)\otimes X^*\oplus S^q
\end{array}.
$$
By Corollary 1, we have an exact sequence
$$
0\ra \oplus S(-d_k)\ra S\otimes Y^*\stackrel{sK^t} {\longrightarrow} S(1)\otimes X^*.
$$
Using this, we can easily see that the kernels of
$$\oplus S(-d_k)  \stackrel{\widetilde{M}^{tr}(0,s)}{\ra}  S^q\ \ \ \tx{and}\ \ \
S\otimes Y^*   \ra S(1)\otimes X^*\oplus S^q $$
 are canonically isomorphic.
 Letting $F_1$ and $F_2$ denote these kernels, from the above commutative diagram, we therefore have the following one
 $$
 \begin{array}{cccc}
 0 \ra &E_1/s_0E_1 & \ra & F_1\\
      &    \da      &     & \da \\
0 \ra & E_2/s_0E_2 & \ra &  F_2
 \end{array},
 $$
 where the downward arrows are isomorphisms. This proves the claim.

The proof of lemma is complete. $\quad\Box$

\section{Minimality}

We begin with the following lemma.

\begin{lem} Let $P$ be a row-homogeneous polynomial matrix of size $p\times q$, and let $a_1,\ldots , a_p$ be nonnegative integers. If
$Q=diag(s_0^{a_1},\ldots ,s_0^{a_p})P$
is regular at infinity, then so is $P$.
\end{lem}
{\em Proof}. Let $f\in T^q$ be such that $s_0f\in ImP^{tr}$, that is,
$$
s_0f=P^{tr}\left(\begin{array}{c}x_1\\x_2\\ \vdots \\x_p\end{array}\right).
$$
We then have
$$
s_0^{1+|a|}f=P^{tr}\left(\begin{array}{c}s_0^{|a|}x_1\\ \vdots \\s_0^{|a|}x_p\end{array}\right)=Q^{tr}\left(\begin{array}{c}s_0^{|a|-a(1)}x_1\\ \vdots \\s_0^{|a|-a(p)}x_p\end{array}\right).
$$
By the hypothesis,
$f\in ImQ^{tr}$. Hence, $f\in ImP^{tr}$.

The proof is complete. $\quad\Box$

Recall that, in \cite{L2}, for every linear differential system ${\cal B}$ and for every AR-model $R$,  have been  introduced important integer functions
$\gamma_{\cal B}$ and $\gamma_R$.

If  ${\cal B}$ is a linear differential system, we set
$$
\mu({\cal B})=\sum_{d\geq 0} \gamma_{\cal B}(d)\frac{nd}{d+1}\cdot \left(\begin{array}{c}
 n+d  \\
d
\end{array}\right)
 \ \  \tx{and}\ \
 \nu({\cal B})=\sum_{d\geq 0} \gamma_{\cal B}(d)\left(\begin{array}{c}
 n+d  \\
d
\end{array}\right).
$$
Likewise, if $R$ is  an AR-model $R$, set
$$
\mu(R)=\sum_{d\geq 0} \gamma_R(d)\frac{nd}{d+1}\cdot \left(\begin{array}{c}
 n+d  \\
d
\end{array}\right)
 \ \  \tx{and}\ \
 \nu(R)=\sum_{d\geq 0} \gamma_R(d)\left(\begin{array}{c}
 n+d  \\
d
\end{array}\right).
$$

\begin{thm} Let $\cal B$  be a linear differential system, and let $\Sigma: \ X\stackrel{K,L}{\longrightarrow} Y \stackrel{M}{\la} \mathbb{F}^q$ be any its
proper KW-representation. Then,
$$
 \mu({\cal B})\leq \dim(X)\ \ \ \tx{and}\ \ \  \nu({\cal B})\leq \dim(Y).
$$
Moreover,
$$
 \mu({\cal B}) = \dim(X)\ \Longleftrightarrow \  \nu({\cal B}) = \dim(Y).
$$
\end{thm}

{\em Proof}. Let $d_1,\ldots ,d_p$ be the Kronecker indices of $\Sigma$, and let $\widetilde{M}$ be  as in Section 5, i.e., the composition $T^q \stackrel{M}{\ra} Y\otimes T \ra \oplus T(d_k)$. Then,
$R=\widetilde{M}^{dh}$ is an
 AR-representation of $\cal B$ (see the proof of Theorem 1). The row degrees of $R^h$ are less than or equal to those of $\widetilde{M}$. Because the latter
  is regular at infinity,
by the previous lemma, $R^h$ also is regular at infinity. Hence, $R$ is a proper AR-representation. Since $\gamma_{\cal B}\leq \gamma_R$, it is clear that
$$
 \mu({\cal B})\leq \mu(R)\ \ \ \tx{and}\ \ \  \nu({\cal B})\leq \nu(R).
$$
Also, it is clear that
$$
 \mu({\cal B}) = \mu(R)\ \Leftrightarrow \ \gamma_{\cal B} = \gamma_R\  \Leftrightarrow \ \nu({\cal B}) = \nu(R).
$$

Further, let $(a_1,\ldots ,a_p)$ be the row degree of $R$. Then $a_i\leq d_i$ for each $i\in [1, p]$, and therefore, in view of (6),
$$
 \mu(R) \leq  \dim(X)\ \ \ \tx{and}\ \ \  \nu(R)\leq \dim(Y).
$$
It is easily seen that
$$
 \mu(R) = \dim(X)\ \Leftrightarrow \ \forall i, \ a_i= d_i \ \Leftrightarrow \ \nu(R) = \dim(Y).
$$

One easily completes the proof. $\quad\Box$

Two proper KW-models are said to be equivalent if they  generate the same manifest behavior.
A KW-model is said to be minimal if among all proper KW-models in its equivalence class it has  minimal number
 of state variables and minimal number of auxiliary variables.

{\em Remark}. The above notion of minimality is adapted from Section VII in Willems \cite{W}.

\begin{thm} Let $\cal B$  be a linear differential system. There exists  a proper KW-representation $X\ra Y \la \mathbb{F}^q$
that satisfies the
 following equivalent conditions:

(a) $\mu({\cal B}) = \dim(X)$;

(b) $\nu({\cal B}) = \dim(Y)$.\\
 Moreover,
such a representation  is uniquely determined up to similarity.
\end{thm}
{\em Proof}. We have already seen that the conditions are equivalent.

"Existence"\ Let $R$ be a minimal proper AR-representation of $\cal B$. Because $\gamma_R=\gamma_{\cal B}$, obviously
$\mu(R)=\mu({\cal B})$ and $\nu(R)=\nu({\cal B})$. Further, because the Kronecker indices of
 $KW(R)$ coincide with  the row degrees of $R$, obviously the dimension of state variables in $KW(R)$ is equal to $\mu(R)$ and the dimension of
  auxiliary variables is equal to $\nu(R)$. It follows that $KW(R)$ satisfies the conditions of theorem.

 "Uniqueness"\  Assume that $\Sigma$ is a minimal proper representation of $\cal B$. Construct an AR-representation $R$ of $\cal B$ as in
 the proof of Theorem 1. It is clear from the proof of Theorem 2 that $\gamma_R=\gamma_{\cal B}$, and hence $R$  is a minimal proper AR-representation.
 Clearly, $\Sigma \simeq KW(R)$.

  Recall now that, by the main result of \cite{L3}, any two minimal proper AR-representations of a linear differential system are Brunovsky equivalent.
  (The reader is referred to the just mentioned paper for the definition of the group $\Gamma(d)$ and the concept of Brunovsky equivalence.)

   We therefore have to show that if $R$ and $R^\prime$ are two Brunovsky equivalent AR-models, then $KW(R)\simeq KW(R^\prime)$. After permuting the
  rows (if necessary),
 we may assume that both of $R$ and $R^\prime$  have the same row degree $d=(d_1,\ldots , d_p)$. Let $U\in \Gamma(d)$ be such that $R^\prime=UR$.
First of all, we  have a commutative diagram
 \begin{equation}
\begin{array}{ccc}
  \mathbb{F}^q & \stackrel{R}{\ra} & \oplus S_{\leq d_i}  \\
 || &  & \da \ U \\
 \mathbb{F}^q & \stackrel{R^\prime}{\ra} & \oplus S_{\leq d_i}
\end{array}.
 \end{equation}
Since $U$ preserves the filtration, for each $k$, the linear map
$$
\oplus S_{\leq d_i+k}\stackrel{U}{\ra} \oplus S_{\leq d_i+k}
$$
is bijective, and if $\lambda$ is a  polynomial of degree $\leq k$ the  diagram
 $$
\begin{array}{ccc}
  \oplus S_{\leq d_i} &\stackrel{\lambda}{\ra} & \oplus S_{\leq d_i+k}  \\
U \da \ \ \ \          &     & \da \ U \\
 \oplus S_{\leq d_i} & \stackrel{\lambda}{\ra} & \oplus S_{\leq d_i+k}
\end{array}
$$
is commutative. Hence, the diagram
 $$
\begin{array}{ccc}
  (\oplus S_{\leq d_i})^{n+1} & \stackrel{\left[\begin{array}{cccc}-1&s_1& \ldots & s_n\end{array}\right]}{\longrightarrow} & \oplus S_{\leq d_i+1}  \\
 & &   \\
 \da                          &     & \da  \\
 (\oplus S_{\leq d_i})^{n+1} & \stackrel{\left[\begin{array}{cccc}-1&s_1& \ldots & s_n\end{array}\right]}{\longrightarrow} & \oplus S_{\leq d_i+1}
\end{array}
$$
is commutative. In view of the exact sequence (4), it follows  that there is a bijective map $\oplus X^{(d_i)}\ra  \oplus X^{(d_i)}$
for which the diagram
$$
\begin{array}{ccccccc}
0\ra & \oplus X^{(d_i)} & \ra & \oplus (S_{\leq d_i})^{n+1} & \ra & \oplus S_{\leq d_i+1}& \ra 0\\
     & \da &     &     \da        &     & \da  & \\
0\ra & \oplus X^{(d_i)} & \ra & \oplus (S_{\leq d_i})^{n+1} & \ra & \oplus S_{\leq d_i+1}& \ra 0\\
\end{array}
$$
commutes.
This together with (8) implies  that $KW(R)\simeq KW(R^\prime)$.

The proof is complete.
  $\quad\Box$

\begin{cor} Every minimal proper KW-model is obtained by applying the functor $KW$ to a minimal proper AR-model.
\end{cor}

\begin{cor} There is a one-to-one correspondence between linear differential systems and similarity classes of minimal proper KW-models.
\end{cor}

\begin{cor} The Kronecker indices of the minimal proper  KW-representation of a linear differential system coincide with the observability indices
of the latter.
\end{cor}

\section{Concluding remark}

In \cite{L3}, we assigned canonically a KW-model to every polynomial matrix, and we have seen here
that this assignment preserves properness and minimality. As a consequence, we get that it
induces a one-to-one correspondence between linear differential systems and similarity
classes of minimal proper KW-models.

\end{document}